\documentclass[11pt]{article}
\usepackage{amsmath,amssymb,theorem}
\usepackage{graphicx, enumerate}
\usepackage{subfigure}
\usepackage{psfrag}
\usepackage{placeins}
\usepackage{wrapfig}
\usepackage{booktabs}
\usepackage{chngcntr}
\counterwithin{figure}{section}
\numberwithin{figure}{section}
\counterwithin{table}{section}
\usepackage{hyperref}
\hypersetup{
  colorlinks,
  citecolor=black,
  filecolor=black,
  linkcolor=black,
  urlcolor=black
}

\newtheorem{thm}{Theorem}[section]
\newtheorem{claim}{Claim}[section]
\newtheorem{conj}[thm]{Conjecture}
\newtheorem{cor}[thm]{Corollary}
\theorembodyfont{\rmfamily}
\newtheorem{rem}[thm]{Remark}
\def\pf{\bigskip\noindent {\bf Proof.}~~}

\def\dfn#1{{\sl #1}}
\def\es{\emptyset}
\def\less{\backslash}

\def\pf{\bigskip\noindent {\emph{Proof.}}~~}

\def\qed{ \hfill $\blacksquare$}

\begin{document}
\title{A note on Hadwiger's Conjecture for $W_5$-free graphs with independence number two}
\author{
Christian Bosse\thanks{Email: csbosse@knights.ucf.edu}
\\
Department  of Mathematics\\
 University of Central Florida\\
Orlando, FL 32816, USA\\
}
\maketitle
\begin{abstract}
The \emph{Hadwiger number} of a graph $G$, denoted $h(G)$, is the largest integer $t$ such that $G$ contains $K_t$ as a minor.  A famous conjecture due to Hadwiger in 1943 states that for every graph $G$, $h(G) \ge \chi(G)$, where $\chi(G)$ denotes the chromatic number of $G$.  Let $\alpha(G)$ denote the independence number of $G$.  A graph is \emph{$H$-free} if it does not contain the graph $H$ as an induced subgraph.  In 2003, Plummer, Stiebitz and Toft proved that $h(G) \ge \chi(G)$ for all $H$-free graphs $G$ with $\alpha(G) \le 2$, where $H$ is any graph on four vertices with $\alpha(H) \le 2$, $H=C_5$, or $H$ is a particular graph on seven vertices.  In 2010, Kriesell considered a particular strengthening of Hadwiger's conjecture due to Seymour and subsequently generalized the statement to include all forbidden subgraphs $H$ on five vertices with $\alpha(H) \le 2$.  In this note, we prove that $h(G) \ge \chi(G)$ for all $W_5$-free graphs $G$ with $\alpha(G) \le 2$, where $W_5$ denotes the wheel on six vertices.  
\end{abstract}

\baselineskip 16pt
\section{Introduction}

All graphs in this note are finite and simple; that is, they have no loops or parallel edges. Given a graph $G$ and a set $X\subseteq V(G)$,  we use $|G|$ to denote  the  number of vertices    of $G$, and  $G[X]$ to denote the subgraph of $G$ obtained from $G$ by deleting all vertices in $V(G)\less X$.  A graph $H$ is an \dfn{induced subgraph} of $G$ if $H=G[X]$ for some $X\subseteq V(G)$.  
For any positive integer $t$, we write  $[t]$ to denote the set $\{1,2, \ldots, t\}$. We use the convention ``$X:=$'' to mean that $X$ is defined to be the right-hand side of the relation.  Given $X, Y \subseteq V(G)$, we say that $X$ is \emph{complete} (resp. \emph{anticomplete}) to $Y$ if for every $x \in X$ and every $y \in Y$, $xy \in E(G)$ (resp. $xy \not \in E(G)$).   A \emph{clique} is a set of pairwise adjacent vertices in $G$, and a set of pairwise non-adjacent vertices is \emph{independent}.  We use $\chi(G)$, $\omega(G)$ and $\alpha(G)$ to denote the chromatic, clique and independence number of a graph $G$, respectively.  An edge $xy \in E(G)$ is \emph{dominating} if every vertex in $G \less \{x, y\}$ is adjacent to either $x$ or $y$. We say that $G$ is \emph{$H$-free} if $G$ contains no induced subgraph isomorphic to the graph $H$.  A graph $H$ is a \emph{minor} of $G$ if $H$ can be obtained from a subgraph of $G$ upon contracting edges.  Define the \emph{Hadwiger number} of a graph $G$, denoted $h(G)$, to be the largest integer $t$ such that $G$ contains a $K_t$ minor.

In 1943 Hadwiger \cite{Hadwiger} made the following famous conjecture.

\begin{conj}\label{HC}
Let $G$ be any graph.  Then $h(G) \ge \chi(G)$.
\end{conj}

Many regard Hadwiger's conjecture as perhaps one of the most profound unsolved problems in graph theory due to its connection with the Four Color Theorem (see \cite{Seymoursurvey} and \cite{Toftsurvey} for more details and history surrounding the conjecture).  To date, a general proof of this conjecture remains elusive.  Several weakenings and special cases have been considered and many partial results lending further credence to it have been obtained.  We refer the reader to Seymour's recent survey \cite{Seymoursurvey} for a fairly complete listing of the partial results.  We simply add here that since the time that survey was published, some new partial results have been obtained.  In 2017, Song and Thomas \cite{ThomasSong} showed that if $\alpha(G) \ge 3$ and $G$ is $\{C_4, C_5, C_6, \ldots, C_{2\alpha(G) - 1}\}$-free, then $h(G) \ge \chi(G)$. Note that Hadwiger's conjecture can be equivalently formulated in the following manner.  For all $t \ge 0$, then $\chi(G) \le t$ for every $K_{t+1}$ minor-free graph.  By $K_t^-$ (resp. $K_t^{=}$), we denote the complete graph $K_t$ with one edge (resp. two edges) removed.  Rolek and Song \cite{RolekSong} showed in 2017 that $\chi(G) \le 8$, $9$ and $12$ for every $K_8^{=}$, $K_8^-$ and $K_9$ minor-free graph, respectively.  Rolek \cite{Rolek} showed later in 2018 that $\chi(G) \le 10$ for every $K_9^{=}$ minor-free graph.

In this note we pay particular attention to Conjecture~\ref{HC} for graphs $G$ with $\alpha(G) \le 2$.  Observe that the complement of any such graph is triangle-free.  As Plummer, Stiebitz and Toft point out in \cite{PST}, this is a mild restriction considering the wide variety of triangle-free graphs.  Seymour says in \cite{Seymoursurvey} the following about Hadwiger's conjecture in this setting.
\begin{quote}
``This seems to me to be an excellent place to look for a counterexample. My own belief is, if it is true for graphs with stability number two then it is probably true in general, so it would be very nice to decide this case."
\end{quote}
Considering graphs with $\alpha(G)\le 2$ thus appears to be a worthwhile playground to explore in order to gain more insight into Conjecture~\ref{HC}.

We first mention a very useful result of Plummer, Stiebitz and Toft \cite{PST}  that establishes an equivalence of Hadwiger's conjecture in this context.


\begin{thm}[\cite{PST}]\label{HC equiv to WC}
Let $G$ be a graph  with $\alpha(G) = 2$.  Then $h(G) \ge \chi(G)$ if and only if $h(G) \ge \lceil |G|/2\rceil$.
\end{thm}

In the same paper, Plummer, Stiebitz and Toft \cite{PST} proved the following. 

\begin{figure}
\centering
\includegraphics[scale=.75]{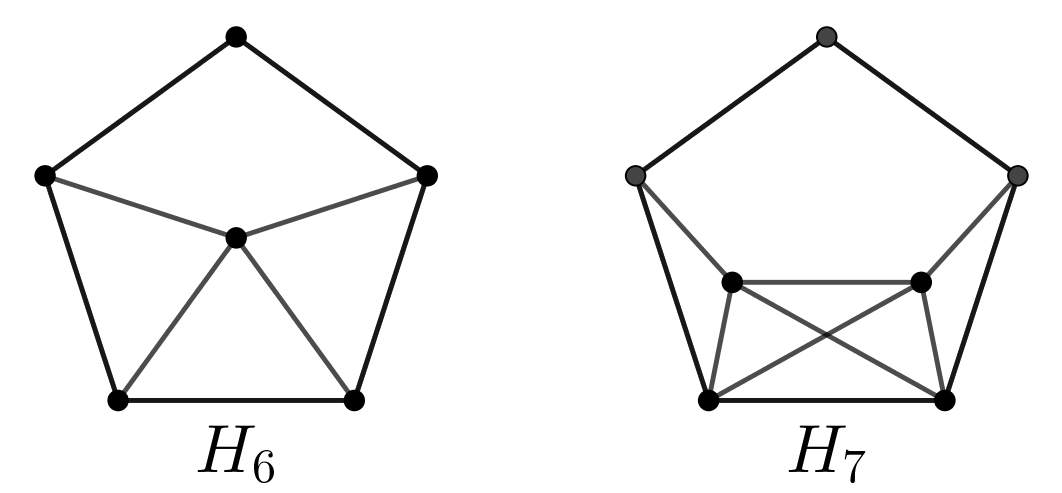}
\caption{The graphs $H_6$ and $H_7$}
\label{H6H7}
\end{figure}

\begin{thm}[\cite{PST}]\label{H7}
Let $G$ be a graph with $\alpha(G) \le 2$.  If $G$ is $H$-free, where   $H$ is  a graph   with  $|H|=4$ and $\alpha(H) \le 2$, or $H=C_5$, or $H=H_7$ (see Figure~\ref{H6H7}),  then $h(G) \ge  \chi(G)$.
\end{thm}

In 2010, Kriesell \cite{Kriesell} further augmented this list of forbidden subgraphs to include all cases of graphs with   five vertices.

\begin{thm}[\cite{Kriesell}]\label{Kriesell}
Let $G$ be a graph with $\alpha(G) \le 2$.  If $G$ is  $H$-free, where   $H$ is  a graph with  $|H|=5$ and $\alpha(H) \le 2$, or $H=H_6$ (see Figure~\ref{H6H7}),  then $h(G) \ge \chi(G)$.
\end{thm}

 Let $W_5:=K_1+C_5$ denote the wheel on six vertices.  In this note, we study Conjecture~\ref{HC} for $W_5$-free graphs with independence number at most two. We prove the following main result. 

\begin{thm}\label{W5}
Let $G$ be a  graph   with $\alpha(G) \le 2$.  If $G$  is  $W_5$-free, then $h(G) \ge \chi(G)$.
\end{thm}

 Our proof of Theorem~\ref{W5} relies only on Theorem~\ref{H7} when $H=C_5$ and   the following result of Chudnovsky and Seymour \cite{seagulls}. 

\begin{thm}[\cite{seagulls}]\label{packingseagull}
Let G be a  graph   with $\alpha(G) \le 2$. If 
\[
\omega(G)\ge \begin{cases}
   |G|/4, & \text{if $|G|$ is even}\\
   (|G|+3)/4, & \text{if $|G|$ is odd,}
\end{cases}
\]
then $h(G) \ge \chi(G)$.
\end{thm}

It is worth noting that if  $G$ is a $K_6$-free graph on $n$ vertices with $\alpha(G) \le 2$  but  does not satisfy Conjecture~\ref{HC},   then $G$ contains a $K_5$ by Theorem~\ref{Kriesell}, and  $n \le 17$ because  $R(K_3, K_6) = 18$ (see \cite{Kery}). But then by Theorem~\ref{packingseagull} and Theorem~\ref{HC equiv to WC},  $h(G) \ge \chi(G)$, a contradiction.  Similarly, if  $G$ is a $K_7$-free graph on $n$ vertices with $\alpha(G) \le 2$  but  does not satisfy Conjecture~\ref{HC},   then  $G$ contains a $K_6$, and $n \le 22$ because   $R(K_3, K_7) = 23$ (see \cite{GraverYackel} and \cite{Kalb}).  But then by Theorem~\ref{packingseagull} and Theorem~\ref{HC equiv to WC},   $h(G) \ge \chi(G)$, a contradiction.  We summarize these observations as follows. 

\begin{rem}
Let $G$ be a $K_t$-free graph with $\alpha(G) \le 2$, where $t \le 7$.  Then $h(G) \ge \chi(G)$.
\end{rem}

We end this section with Corollary~\ref{W5Cor} below, which  follows   from  Theorem~\ref{W5}. 

\begin{cor}\label{W5Cor}
Let $G$ be a graph with $\alpha(G) \le 2$.  If $G$ is $\overline{K_{1,5}}$-free, then $h(G) \ge \chi(G)$.
\end{cor}

\pf Let $G$ be a $\overline{K_{1,5}}$-free graph on $n$ vertices with $\alpha(G)\le2$. By Theorem~\ref{HC equiv to WC}, it suffices to show that  $h(G) \ge  \lceil n/2\rceil$.   Suppose   $h(G) < \lceil n/2\rceil$.  By Theorem~\ref{W5}, $G$ must contain an induced $W_5$, say with vertices    $x_1, x_2, x_3, x_4, x_5, z$, where $G[\{x_1, x_2, x_3, x_4, x_5\}]=C_5$.   We choose such a graph $G$ with $n$ minimum.   By the minimality of $n$, $G$ has no dominating edges.    For all $i \in [5]$, since $zx_i$ is not a dominating edge,  there must exist a vertex  $y_i \in V(G)\less V(W_5)$ such that $y_ix_i, y_iz \not \in E(G)$.  Then $y_ix_{i+2}, y_ix_{i+3} \in E(G)$ because   $\alpha(G) = 2$, where  all arithmetic on indices here   is done
modulo  $5$.    It follows that for all $i, j\in[5]$ with $i\ne j$, $y_i\ne y_j$. Furthermore, $G[\{y_1, y_2, y_3, y_4, y_5\}]=K_5$.  But then $G[\{y_1, y_2, y_3, y_4, y_5, z\}] = \overline{K_{1,5}}$, a contradiction. This completes the proof of Corollary~\ref{W5Cor}.
  \qed\\

 We are now  ready to prove  Theorem~\ref{W5}  in Section~2.
\FloatBarrier

\section{Proof of Theorem~\ref{W5}}
Let $G$ be a $W_5$-free graph on $n$ vertices with $\alpha(G)\le2$. By Theorem~\ref{HC equiv to WC}, it suffices to show that  $h(G) \ge  \lceil n/2\rceil$.   Suppose   $h(G) < \lceil n/2\rceil$.    By Theorem~\ref{H7}, $G$ must contain an induced $C_5$.   We choose such a graph $G$ with $n$ minimum.    Then   $\alpha := \alpha(G)=2$.  Note that $(n+3)/4 \le \lceil (n+2)/4\rceil$ for odd $n$.  By Theorem~\ref{packingseagull},  $\omega(G)< \lceil (n+2)/4 \rceil$ when $n$ is odd, and $\omega(G) < \lceil n/4 \rceil$ when $n$ is even.\medskip

Since $G$ has an induced $C_5$, let $X:= \bigcup_{i=1}^5 X_i$ be a maximal inflation of $C_5$ in $G$ such that   for all $i\in[5]$, $X_i$ is a clique;     $X_i$ is complete to $X_{i-1} \cup X_{i+1}$, and anticomplete to $X_{i-2} \cup X_{i+2}$, where  all arithmetic on indices here and henceforth is done
modulo  $5$.  Then  $X_i \ne \es$ for all $i \in [5]$  and  $G[X_i]$ is a clique for every $i \in [5]$.  Since $\alpha=2$ and $G$ is $W_5$-free, no vertex in $G$ is complete to $X$ and every vertex in $G  \less X$ must be complete to at least three consecutive $X_i$'s on the maximal inflation of $C_5$. 
For each $i\in[5]$, let 
 \[
\begin{split}
Y_i  &: =\{v \in V(G) \less X \mid  v \text{ is complete to } X\less X_i \text{ and has a non-neighbor in } X_i \}\\
Z_i  &: =\{v \in V(G) \less X \mid  v \text{ is complete to } X\less (X_i\cup X_{i+1}) \text{ and has a non-neighbor in } X_i \text{ and   in } X_{i+1}\}  
\end{split}
\]
Let  $Y:=\bigcup_{i=1}^5 Y_i$ and  $Z:=\bigcup_{i=1}^5 Z_i$.  By definition, $Y \cap Z = \es$ and $Y\cup Z=V(G)\less X$. By the maximality of $|X|$, no vertex in $Z_i$ is anticomplete to $X_i\cup X_{i+1}$ in $G$, else, such a vertex can be placed in $X_{i+3}$ to obtain a larger inflation of $C_5$.  

\begin{claim}\label{Y_i cliques}
For all $i \in [5]$, $Y_i$ is anticomplete to $X_i$, and so  $G[Y_i]$ is a clique.
\end{claim}
\begin{pf}
Suppose some $Y_i$, say $Y_1$ is not anticomplete to $X_1$.  Then $Y_1\ne\es$.   Since $Y_1$ is not anticomplete to $X_1$, there exist   $y_1 \in Y_1$ and $x_1, x_1'\in X_1$  such that $y_1x_1 \not\in E(G)$ but   $y_1x_1' \in E(G)$.  Let $x_i \in X_i$ for all $i \in \{2, 3, 4, 5\}$.  By definition, $y_1x_i \in E(G)$ for all $i \in \{2, 3, 4, 5\}$.  But then $G[\{y_1, x_1', x_2, x_3, x_4, x_5\}]=W_5$, a contradiction.   Thus for all $i\in[5]$, $Y_i$ must be  anticomplete to $X_i$.  It follows that  $G[Y_i]$ must be a clique because $\alpha = 2$.\qed
\end{pf}

\begin{claim}
For all $i \in [5]$, $G[Z_i]$ is a clique.
\end{claim}
\begin{pf}
Suppose some $G[Z_i]$, say $G[Z_1]$, is not a clique.    Then there exist   $z_1, z_1' \in Z_1$ such that $z_1z_1' \not\in E(G)$.  By definition of $Z_1$, there exist $x_1 \in X_1$ and $x_2 \in X_2$ such that  $z_1x_1, z_1x_2 \not\in E(G)$. Since $\alpha = 2$, we see that  $z_1'x_1, z_1'x_2 \in E(G)$.  But then $G[\{z_1', x_1, x_2, x_3, x_4, x_5\}] =W_5$, where $x_i \in X_i$ for all $i \in \{3, 4, 5\}$, a contradiction. \qed 
\end{pf}

\begin{claim}\label{Z_i}
For all $i\in[5]$, every vertex in $  Z_i$ is  either   anticomplete to $X_i$,  or   anticomplete to $X_{i+1}$, but not both.
\end{claim}
\begin{pf}
As observed earlier, for all $i\in[5]$, no vertex in $Z_i$ is anticomplete to $X_i\cup X_{i+1}$.  Suppose there exists some $i\in[5]$, say $i=1$, such that some vertex, say $z\in Z_1$ is neither  anticomplete to $X_i$   nor   anticomplete to $X_{i+1}$.  Then there exist $x_1, x_1' \in X_1$ and $x_2, x_2' \in X_2$ such that 
  $zx_1, zx_2\not\in E(G)$ and $zx_1', zx_2' \in E(G)$.  Let $x_i \in X_i$ for all $i \in \{3, 4, 5\}$. By definition of $Z_1$, $z$ is complete to $\{x_3, x_4, x_5\}$. But then   $G[\{z, x_1', x_2', x_3, x_4, x_5\}]=W_5$,   a contradiction.\qed \\
\end{pf}

For each $i\in[5]$, let 
 \[
\begin{split}
Z_i^i  &: =\{z \in Z_i \mid  z \text{ is anticomplete to }  X_i   \}\\
Z_i^{i+1}  &: =\{z \in Z_i \mid  z \text{ is anticomplete to }  X_{i+1}\}  
\end{split}
\]
By Claim~\ref{Z_i},  $Z_i = Z_i^i \cup Z_i^{i+1}$ and $Z_i^i \cap Z_i^{i+1}=\es$ for all $i\in[5]$. \\

 Since   $\alpha=2$,  by the choice of $Y_i, Z_i, Z_i^i, Z_i^{i+1}$, we see   that 

\begin{claim}\label{ZYZ cliques}For all $i \in [5]$, 
both $G[Z_{i-1}^i \cup Y_i \cup Z_i]$ and $G[Z_{i-1} \cup Y_i \cup Z_i^i]$  are cliques. 
\end{claim}
 
 We next show that 
 

\begin{claim}\label{Zi^i parts}
For all $i \in [5]$, every vertex in $Z_i^i$ is complete to $Y_{i-1}$ or complete to $Z_{i+1}^{i+2}$.
\end{claim}
\begin{pf}  
Suppose the statement is false.  We may assume that there exists  some vertex $z \in Z_1^1$ such that $zy_5, zz_2 \not \in E(G)$, where $y_5 \in Y_5$ and $z_2 \in Z_2^3$.  Since $\alpha = 2$, we see that  $y_5z_2 \in E(G)$.  Then $G[\{y_5, z_2, x_5, z, x_3, x_4\}] = W_5$, where $x_5\in X_5$, $x_3\in X_3$, $x_4\in X_4$, a contradiction. \qed  
\end{pf}

\begin{claim}\label{Y_i partition}
For all $i \in [5]$,    every vertex in  $Y_i$ is either  complete to $Y_{i-1}$ or complete to $Y_{i+2}$.
\end{claim}
\begin{pf}
Suppose not.  We may assume there exist  vertices $y_1 \in Y_1$, $y_3 \in Y_3$ and $y_5 \in Y_5$ such that $y_1y_3, y_1y_5 \not\in E(G)$.  Then $y_3y_5 \in E(G)$ because $\alpha = 2$.  Then $G[\{y_5, y_3, x_5, y_1, x_3, x_4\}] = W_5$, where   $x_5\in X_5, x_3\in X_3, x_4\in X_4$, a contradiction. \qed\\
\end{pf}

By Claim~\ref{Zi^i parts},   $Z_1^1 = A_1 \cup B_1$ and $Z_3^3 = A_3 \cup B_3$, where $A_i \cap B_i = \es$, and 
\begin{align*}
A_i &:= \{v \in Z_i^i \mid \text{$v$ is complete to $Y_{i-1}$}\}\\
B_i &:= \{v \in Z_i^i \mid \text{$v$ is complete to $Z_{i+1}^{i+2}$ and has a non-neighbor in  $Y_{i-1}$}\}.
\end{align*}
for $i \in \{1, 3\}$.    By Claim~\ref{Y_i partition},   $Y_1 = Y_1' \cup Y_1''$, where
\begin{align*}
Y_1' &:= \{v \in Y_1 \mid \text{$v$ is complete to $Y_5$}\}\\
Y_1'' &:= \{v \in Y_1 \mid \text{$v$ is complete to $Y_3$ and has a non-neighbor  in $Y_5$}\}.
\end{align*}
Then  $Y_1' \cap Y_1'' = \es$.  We claim that $A_3$ is complete to $Y_1''$ in $G$.  To see this, suppose there exist vertices $z \in A_3$ and $y_1 \in Y_1''$ such that $zy_1 \not\in E(G)$.  
By the choice of $Y_1''$, there exists a vertex $y_5\in Y_5$ such that  $y_1y_5 \not\in E(G)$. Then  $zy_5 \in E(G)$ because $\alpha = 2$.  Since $z \in Z_3^3$, there exists some vertex $x_4 \in X_4$ such that $zx_4 \in E(G)$.  But then $G[\{z, y_5, x_3, y_1, x_5, x_4\}] = W_5$, where $x_3\in X_3$ and $x_5\in X_5$, a contradiction. This proves that $A_3$ is complete to $Y_1''$ in $G$, as claimed.  Let 
\begin{align*}
H_1 &:= G[X_3 \cup X_4 \cup Y_5 \cup Z_5 \cup Y_1' \cup A_1]\\
H_2 &:= G[X_4 \cup X_5 \cup B_1 \cup Z_1^2 \cup Y_2 \cup Z_2]\\
H_3 &:= G[X_1 \cup X_2 \cup B_3 \cup Z_3^4 \cup Y_4 \cup Z_4]\\
H_4 &:= G[X_5 \cup Y_1'' \cup Y_3 \cup A_3]
\end{align*}
Note that each of $H_1$, $H_2$, $H_3$ and $H_4$ is a clique in $G$, and  $|H_1| + H_2| + |H_3| +|H_4|=|G|+|X_4|+|X_5|\ge n+2$. It follows that $\omega(G)\ge \max\{|H_1|, |H_2|, |H_3|, |H_4|\}\ge \lceil (n+2)/4\rceil  $, a contradiction.\medskip 

This completes the proof of Theorem~\ref{W5}. \qed

\section{Acknowledgement}
I would like to thank my advisor Zi-Xia Song for her help and guidance on this topic.

\bibliographystyle{amsplain}
\bibliography{refs}

%
%
%
%
\end{document}